\documentclass[12pt,a4paper,english,intoc,bibliography=totoc,index=totoc,BCOR10mm,captions=tableheading,titlepage,fleqn]{scrbook}
\usepackage{lmodern}

\usepackage[T1]{fontenc}
\usepackage[latin9]{inputenc}
\usepackage{babel}
\usepackage{amsmath}
\usepackage{amssymb}
\usepackage[unicode=true,
 bookmarks=true,bookmarksnumbered=true,bookmarksopen=true,bookmarksopenlevel=2,
 breaklinks=false,pdfborder={0 0 0},backref=false,colorlinks=false]
 {hyperref}
\hypersetup{pdftitle={Abstract},
 pdfauthor={Uwe Stöhr},
 pdfpagelayout=OneColumn, pdfnewwindow=true, pdfstartview=XYZ, plainpages=false}
\usepackage{breakurl}

\makeatletter

\usepackage{calc}

\makeatother

\begin{document}
\noindent \begin{center}
\textbf{\Large Mad families and non-meager filters}
\par\end{center}{\Large \par}

\noindent \begin{center}
{\large Haim Horowitz and Saharon Shelah}
\par\end{center}{\large \par}

\noindent \begin{center}
\textbf{Abstract}
\par\end{center}

\noindent \begin{center}
{\small We prove the consistency of $ZF+DC+"$there are no mad families$"+"$there
exists a non-meager filter on $\omega"$ relative to $ZFC$, answering
a question of Neeman and Norwood. We also introduce a weaker version
of madness, and we strengthen the result from {[}HwSh:1090{]} by showing
that no such families exist in our model.}%
\footnote{{\small Date: January 11, 2017}{\small \par}



Publication 1103 of the second author

Partially supported by European Research Council grant 338821.%
}
\par\end{center}{\small \par}

\noindent \begin{center}
\textbf{Introduction}
\par\end{center}

This paper is a continuation of {[}HwSh:1090{]}, which is part of
the ongoing effort to investigate the possible non-existence and definability
of mad families. In {[}HwSh:1090{]} we proved that $ZF+DC+"$there
are no mad families$"$ is equiconsistent with $ZFC$ (previous results
by Mathias and Toernquist established the consistency of that statement
relative to large cardinals, see {[}Ma1{]} and {[}To{]}). In this
paper we extend our results from {[}HwSh:1090{]} to address the following
question by Neeman and Norwood:

\textbf{Question ({[}NN{]}): }If there are no mad families, does it
follow that every filter is meager?

By a result of Mathias ({[}Ma2{]}), if every set of reals has the
Ramsey property, then every filter is meager.

We shall construct a model of $ZF+DC$ where there are no mad families,
but there is a non-meager filter on $\omega$. Our proof relies heavily
on {[}HwSh:1090{]}, the main change is that now we're dealing with
a class $K_2$ consisting of pairs $(\mathbb P,\mathcal{A})$ such
that $\mathbb P$ is ccc and forces $MA_{\aleph_1}$, and in addition,
$\mathbb P$ forces that $\mathcal A$ is independent (we shall require
more, see definition 2). In order to imitate the proof from {[}HwSh:1090{]},
we need to prove analogous amalgamation results for an appropriate
subclass of $K_2$. As in {[}HwSh:1090{]}, our final model is obtained
by forcing with $\mathbb P$ where $(\mathbb P,\mathcal A)$ is a
{}``very large'' object in a subclass of $K_2$, and the non-meager
filter will be constructed from $\mathcal A$, which should contain
many Cohen reals.

Finally, we consider the notion of nearly mad families (see definition
14), which was also introduced in {[}NN{]}. We introduce the notion
of a somewhat mad family, which includes both mad and nearly mad families,
and we prove that no somewhat mad families exist in our model. 

\noindent \begin{center}
\textbf{A non-meager filter without mad families}
\par\end{center}

\textbf{Hypothesis 1: }We fix $\mu$ and $\lambda$ such that $\aleph_2 \leq \mu$,
$\lambda=\lambda^{<\mu}$, $\mu=cf(\mu)$ and $\alpha<\mu \rightarrow |\alpha|^{\aleph_1}<\mu$.

\textbf{Definition 2: }A.\textbf{ }Let $K_2$ be the class of $\bold k$
such that:

a. Each $\bold k$ has the form $(\mathbb P,\mathcal A)=(\mathbb{P}_{\bold k}, \mathcal{A_{\bold k}})$.

b. $\mathbb P$ is a ccc forcing such that $\Vdash_{\mathbb P} MA_{\aleph_1}$.

c. $\mathcal A$ is a set of canonical $\mathbb P$-names of subsets
of $\omega$.

d. $\Vdash_{\mathbb P} "\mathcal A$ is independent, i.e. every finite
non-trivial Boolean combination of elements of $\mathcal A$ is infinite$"$.

B. For $\bold k \in K_2$ and a $\mathbb{P}_{\bold k}$-name $\underset{\sim}{b}$,
let $\Vdash_{\mathbb{P}_{\bold k}} "\underset{\sim}{b} \in pos(\bold k)"$
mean $\Vdash_{\mathbb{P}_{\bold k}} "\underset{\sim}{b} \in [\omega]^{\omega}$and
there is no non-trivial Boolean combination of sets from $\mathcal{A}_{\bold k}$
that is almost disjoint to $\underset{\sim}{b}"$.

C. Let $\leq_1$ be the following partial order on $K_2$: $\bold{k}_1 \leq_1 \bold{k}_2$
if and only if:

a. $\mathbb{P}_{\bold{k}_1} \lessdot \mathbb{P}_{\bold{k}_2}$.

b. $\mathcal{A}_{\bold{k}_1} \subseteq \mathcal{A}_{\bold{k}_2}$.

D. Let $\leq_2$ be the following partial order on $K_2$:

$\bold{k}_1 \leq_2 \bold{k}_2$ if and only if:

a. As in B(a).

b. As in B(b).

c. If $\underset{\sim}{b}$ is a $\mathbb{P}_{\bold{k}_1}$-name then
$\Vdash_{\mathbb{P}_{\bold{k}_1}} "\underset{\sim}{b} \in pos(\mathcal{A}_{\bold{k}_1})"$
implies $\Vdash_{\mathbb{P}_{\bold{k}_2}} "\underset{\sim}{b} \in pos(\mathcal{A}_{\bold{k}_2})"$.

E. Let $K_2^+$ be the class of $\bold k \in K_2$ such that $\Vdash_{\mathbb{P}_{\bold k}} "\mathcal{A}_{\bold k}$
is a maximal independent set everywhere$"$, where $\mathcal A$ is
a maximal independent set everywhere if for every $a_0,...,a_{n-1} \in \mathcal A$
without repetition, $b:=\underset{l<n}{\cap}a_l^{\text{if l is even}} \in [\omega]^{\omega}$
and $\{a\cap b : a\in \mathcal A \setminus \{a_0,...,a_{n-1}\}\}$
is a maximal independent set in $[b]^{\omega}$.

F. When we write $"\underset{\sim}{a_0},...,\underset{\sim}{a_{n-1}} \in \mathcal A"$,
we mean that $(\underset{\sim}{a_i} : i<n)$ is without repetition,
moreover, $i<j<n \rightarrow \Vdash_{\mathbb P} "\underset{\sim}{a_i} \neq \underset{\sim}{a_j}$.

\textbf{Observation 3: }a.\textbf{ }$\leq_1$ and $\leq_2$ are partial
orders, and if $\bold{k}_1,\bold{k}_2 \in K_2^+$ then $\bold{k}_1 \leq_1 \bold{k}_2 \rightarrow \bold{k}_1 \leq_2 \bold{k}_2$.

b. If $\bold{k}_1 \leq_2 \bold{k}_2$ and $\underset{\sim}{b}$ is
a $\mathbb{P}_{\bold{k}_1}$-name, then $\Vdash_{\mathbb{P}_{\bold{k}_1}} "\underset{\sim}{b} \in pos(\mathcal{A}_{\bold{k}_1})"$
iff $\Vdash_{\mathbb{P}_{\bold{k}_2}} "\underset{\sim}{b} \in pos(\mathcal{A}_{\bold{k}_2})"$.

\textbf{Proof: }We shall prove the second claim of 3(a), everything
else should be clear. Suppose that $\Vdash_{\mathbb{P}_{\bold{k}_1}} "\underset{\sim}{b} \in pos(\mathcal{A}_{\bold{k}_1})"$,
but for some $\underset{\sim}{a_0},...,\underset{\sim}{a_{n-1}} \in \mathcal{A}_{\bold{k}_2}$
and $p\in \mathbb{P}_{\bold{k}_2}$, $p\Vdash_{\mathbb{P}_{\bold{k}_2}} "\underset{\sim}{b} \cap (\underset{l<n}{\cap} \underset{\sim}{a_l^{\text{l is even}}})$
is finite$"$. Let $G\subseteq \mathbb{P}_{\bold{k}_1}$ be generic
over $V$ such that $p\in G$ and we shall work over $V[G]$. WLOG
there is $n_1<n$ such that $\underset{\sim}{a_l} \in \mathcal{A}_{\bold{k}_1}$
iff $l<n_1$, and denote $\underset{\sim}{a_*}=\underset{l<n_1}{\cap} \underset{\sim}{a_l^{\text{l is even}}}$.
As $\Vdash_{\mathbb{P}_{\bold{k}_1}} "\underset{\sim}{b} \in pos(\mathcal{A}_{\bold{k}_1})"$,
it follows that $\underset{\sim}{b} \cap \underset{\sim}{a_*}$ is
infinite. It's enough to show that for some Boolean combination $\underset{\sim}{a_{**}}$
from $\mathcal{A}_{\bold{k}_1}$, $\underset{\sim}{a_{**}} \subseteq \underset{\sim}{a_*}$
and $\underset{\sim}{a_{**}} \subseteq^* \underset{\sim}{b}$, as
then $\underset{\sim}{a_{**}} \cap (\underset{n_1 \leq l<n}{\cap} \underset{\sim}{a_l^{\text{l is even}}}) \subseteq^* \underset{\sim}{b} \cap (\underset{l<n}{\cap} \underset{\sim}{a_l^{\text{l is even}}})$,
and therefore it's finite, contradicting the definition of $\mathcal{A}_{\bold{k}_2}$.
As $\bold{k}_1 \in K_2^+$, \inputencoding{latin1}{it follows that
}\inputencoding{latin9}$\{ \underset{\sim}{a_*} \cap \underset{\sim}{c} : \underset{\sim}{c} \in \mathcal{A}_{\bold{k}_1} \setminus \{ \underset{\sim}{a_l} : l<n_1\}\}$
is a maximal independent set in $[\underset{\sim}{a_*}]^{\omega}$,
hence there are $\underset{\sim}{c_0},...,\underset{\sim}{c_{m-1}} \in \mathcal{A}_{\bold{k}_1} \setminus \{\underset{\sim}{a_l} : l<n_1\}$
such that $(\underset{l<n_1}{\cap} \underset{\sim}{a_l^{\text{if l is even}}}) \cap (\underset{k<m}{\cap} \underset{\sim}{c_m^{\text{if m is even}}}) \subseteq^* \underset{\sim}{b}$,
so $\underset{\sim}{a_{**}}=(\underset{l<n_1}{\cap} \underset{\sim}{a_l^{\text{if l is even}}}) \cap (\underset{k<m}{\cap} \underset{\sim}{c_m^{\text{if m is even}}})$
is as required. $\square$

\textbf{Observation 4: }$\bold{k}_1 \leq_2 \bold{k}_2$ and $\bold{k}_1 \in K_2^+$
when the following hold for some $\kappa$:

a. $\bold{k}_2 \in K_2^+$.

b. $\bold{k}_2 \in H(\kappa)$.

c. $M$ is a model such that $\bold{k}_2 \in M \prec_{\mathcal{L}_{\aleph_2,\aleph_2}} (H(\kappa),\in)$.

d. $\bold{k}_1=\bold{k}_2^M$.

\textbf{Proof: }By observation 3, recalling that $"\mathbb P \models ccc"$
and $"\mathbb P \models MA_{\aleph_1}"$ are $\mathcal{L}_{\aleph_2,\aleph_2}$-expressible
.

\textbf{Claim 5: }For every $\bold{k} \in K_2$ there is $\bold{k}' \in K_2^+$
such that $\bold{k} \leq_1 \bold{k}'$. Moreover, if $|\mathbb{P}_{\bold k}|<\mu$
then we can find such $\bold{k}'$ that satisfies $|\mathbb{P}_{\bold{k}'}|<\mu$.

\textbf{Proof: }If $|\mathbb{P}_{\bold k}|<\mu$, let $\lambda_*=\mu$,
otherwise, let $\lambda_*$ be a regular cardinal greater than $(2+|\mathbb{P}_{\bold{k}}|)^{\aleph_1}$,
such that $\alpha<\lambda_* \rightarrow |\alpha|^{\aleph_1}<\lambda_*$.

we try to choose a sequence $(\bold{k}_{\alpha} : \alpha<\lambda_*)$
by induction on $\alpha<\lambda_*$ such that:

1. $\bold{k}_0=\bold k$.

2. $(\bold{k}_{\beta} : \beta \leq \alpha)$ is an increasing continuous
sequence of members of $K_2$ (with respect to $\leq_1$).

3. $|\mathbb{P}_{\bold{k}_{\alpha}}| <\mu$.

4. For every $\alpha<\lambda_*$, if $\bold{k}_{\alpha} \notin K_2^+$,
we choose a Boolean combination $\underset{\sim}{a_{\alpha}}$ from
$\mathcal{A}_{\bold{k}_{\alpha}}$ and $\underset{\sim}{b_{\alpha}} \subseteq \underset{\sim}{a_{\alpha}}$
witnessing the failure of the condition from Definition 2(E). We then
define $\mathbb{P}_{\bold{k}_{\alpha+1}}$ as an extension (with respect
to $\lessdot$) of $\mathbb{P}_{\bold{k}_{\alpha}} \star Cohen$ to
a ccc forcing that forces $MA_{\aleph_1}$, we let $\underset{\sim}{\eta_{\alpha}}$
be the relevant Cohen generic real and we let $\mathcal{A}_{\bold{k}_{\alpha+1}}=\mathcal{A}_{\bold{k}_{}\alpha} \cup \{\underset{\sim}{b_{\alpha}} \cup \underset{\sim}{\eta_{\alpha}^{-1}}(\{1\})\}$.

5. If $\alpha<\lambda_*$ is a limit ordinal, we define $\bold{k}_{\alpha}$
as in the proof of claim 7 below.

Why can we carry the induction at stage $\alpha+1$ for $\alpha$
as in (4)? We shall prove that for each $\alpha$, $\Vdash_{\mathbb{P}_{\bold{k}_{\alpha+1}}} "\mathcal{A}_{\bold{k}_{\alpha+1}}$
is independent$"$. Let $\underset{\sim}{a_{\alpha}^*}=\underset{\sim}{b_{\alpha}} \cup \underset{\sim}{\eta_{\alpha}^{-1}}(\{1\})$,
note that $\Vdash_{\mathbb{P}_{\bold{k}_{\alpha+1}}} "\underset{\sim}{a_{\alpha}^*} \notin \mathcal{A}_{\bold{k}_{\alpha}}"$,
as otherwise there are $p\in \mathbb{P}_{\bold{k}_{\alpha+1}}$, $n<\omega$
and $\underset{\sim}{a'} \in \mathcal{A}_{\bold{k}_{\alpha}}$ such
that $p\Vdash_{\mathbb{P}_{\bold{k}_{\alpha+1}}} "\underset{\sim}{a'} \setminus n = \underset{\sim}{a_{\alpha}^*} \setminus n"$,
and therefore $p\Vdash_{\mathbb{P}_{\bold{k}_{\alpha+1}}} "\underset{\sim}{\eta_{\alpha}^{-1}}(\{1\}) \restriction (\omega \setminus \underset{\sim}{b_{\alpha}} \setminus n)=\underset{\sim}{a_{\alpha}^*} \setminus \underset{\sim}{b_{\alpha}} \setminus n = \underset{\sim}{a'} \setminus \underset{\sim}{b_{\alpha}} \setminus n \in V^{\mathbb{P}_{\bold{k}_{\alpha}}}"$,
a contradiction (as $\underset{\sim}{\eta_{\alpha}}$ is Cohen and
$\omega \setminus \underset{\sim}{b_{\alpha}} \setminus n$ is infinite).
Now if $\underset{\sim}{a_{\alpha}}=\underset{l<n}{\cap} \underset{\sim}{a_{\alpha,l}^{\text{if l is even}}}$
and $\underset{\sim}{c}=\underset{\sim}{a_{\alpha}} \cap (\underset{l<m}{\cap} \underset{\sim}{d_l^{\text{if l is even}}})$
where $\underset{\sim}{d_0},...,\underset{\sim}{d_{m-1}} \in \mathcal{A}_{\bold{k}_{\alpha}} \setminus \{\underset{\sim}{a_{\alpha,0}},...,\underset{\sim}{a_{\alpha,n-1}}\}$,
then $\underset{\sim}{c} \setminus \underset{\sim}{a_{\alpha}^*}$
and $\underset{\sim}{c} \cap \underset{\sim}{a_{\alpha}^*}$ are infinite
(as $\underset{\sim}{\eta_{\alpha}}$ is Cohen and $\underset{\sim}{c} \setminus \underset{\sim}{b}$
is infinite), so $\mathcal{A}_{\bold{k}_{\alpha+1}}$ is forced to
be independent.

If for some $\alpha<\lambda_*$, $\bold{k}_{\alpha} \in K_2^+$, when
we're done. Otherwise, by Fodor's lemma, there are $\alpha<\beta<\lambda_*$
such that $(\underset{\sim}{a_{\alpha}},\underset{\sim}{b_{\alpha}})=(\underset{\sim}{a_{\beta}},\underset{\sim}{b_{\beta}})$,
a contradiction. $\square$

\textbf{Definition 6: }We say that $(\bold{k}_{\alpha} : \alpha<\beta)$
is increasing continuous if $\alpha_1<\alpha_1 \rightarrow \bold{k}_{\alpha_1} \leq_2 \bold{k}_{\alpha_2}$,
and for every limit $\delta<\beta$, $\underset{i<\delta}{\cup}\mathbb{P}_{\bold{k}_i} \lessdot \mathbb{P}_{\bold{k}_{\delta}}$.

\textbf{Claim 7: }Every increasing continuous sequence in $(K_2^+,\leq_2)$
has an upper bound. Moreover, if the length of the sequence has cofinality
$>\aleph_1$, then the union is an upper bound in $K_2^+$.

\textbf{Proof: }Given an increasing continuous sequence $(\bold{k}_{\alpha} : \alpha<\beta)$,
we choose $\mathbb{P}_{\bold{k}_{\beta}}$ as in {[}HwSh:1090{]} and
we let $\mathcal{A}_{\bold{k}_{\beta}}=\underset{\alpha<\beta}{\cup}\mathcal{A}_{\bold{k}_{\alpha}}$.
This is enough for $\leq_1$, so by claim 5 we're done. $\square$

\textbf{Claim 8: }A. If $\bold{k}_1 \in K_2$ then there are $\bold{k}_2$
and $\underset{\sim}{a}$ such that:

a. $\bold{k}_1 \leq_1 \bold{k}_2$.

b. $\mathcal{A}_{\bold{k}_1} \cup \{\underset{\sim}{a}\} \subseteq \mathcal{A}_{\bold{k}_2}$.

c. $\underset{\sim}{a}$ is Cohen over $V^{\mathbb{P}_{\bold{k}_1}}$.

d. $|\mathbb{P}_{\bold{k}_2}| \leq (2+|\mathbb{P}_{\bold{k}_1}|)^{\aleph_1}$.

B. Moreover, we may require that $\bold{k}_2 \in K_2^+$.

\textbf{Proof: }A. Let $\mathbb{P}=\mathbb{P}_{\bold{k}_1} \star \mathbb C$
where $\mathbb C$ is Cohen forcing, now let $\mathbb{P}_{\bold{k}_2}$
be a ccc forcing such that $\mathbb{P} \lessdot \mathbb{P}_{\bold{k}_2}$,
$\Vdash_{\mathbb{P}_{\bold{k}_2}} "MA_{\aleph_1}"$ and $|\mathbb{P}_{\bold{k}_2}| \leq (2+|\mathbb{P}|)^{\aleph_1}$.
Finally, let $\mathcal{A}_{\bold{k}_2}=\mathcal{A}_{\bold{k}_1} \cup \{\underset{\sim}{a}\}$
where $\underset{\sim}{a}$ is a name for a Cohen real added by $\mathbb{P}_{\bold{k}_2}$,
it's easy to see that $(\mathbb{P}_{\bold{k}_2},\mathcal{A}_{\bold{k}_2})$
are as required.

B. By claim 5. $\square$

\textbf{Definition 9: }We define the amalgamation property in the
context of $K_2^+$ as follows: $K_2^+$ has the amalgamation property
if A implies B where:

A. a. $\bold{k}_l \in K_2^+$ $(l=0,1,2)$.

b. $\bold{k}_0 \leq_2 \bold{k}_l$ $(l=1,2)$.

c. $\mathbb{P}_{\bold{k}_1} \cap \mathbb{P}_{\bold{k}_2}=\mathbb{P}_{\bold{k}_0}$.

B. There exists $\bold{k}_3=(\mathbb{P}_{\bold{k}_3},\mathcal{A}_{\bold{k}_3}) \in K_2^+$
such that $\bold{k}_l \leq_2 \bold{k}_3$ $(l=1,2)$.

\textbf{Claim 10: }a. $(K_2^+,\leq_2)$ has the amalgamation property.

b. Suppose that $\bold{k}_0,\bold{k}_1$ and $\bold{k}_2 \in K_2^+$,
$g : \mathbb{P}_{\bold{k}_0} \rightarrow \mathbb{P}_{\bold{k}_1}$
is an embedding such that $(g(\mathbb{P}_{\bold{k}_0}),g(\mathcal{A}_{\bold{k}_0})) \leq_2 \bold{k}_1$
and $\bold{k}_0 \leq_2 \bold{k}_2$, then there exist $\bold k, \bold{k}_2' \in K_2^+$
and $f$ such that:

1. $\bold{k}_1 \leq_2 \bold{k}$ and $|\mathbb{P}_{\bold k}| \leq (2+|\mathbb{P}_{\bold{k}_1}|+|\mathbb{P}_{\bold{k}_2}|)^{\aleph_1}$.

2. $(g(\mathbb{P}_{\bold{k}_0}),g(\mathcal{A}_{\bold{k}_0})) \leq_2 \bold{k}_2' \leq_2 \bold{k}$.

3. $f: \mathbb{P}_{\bold{k}_2} \rightarrow \mathbb{P}_{\bold{k}_2'}$
is an isomorphism mapping $\mathcal{A}_{\bold{k}_2}$ to $\mathcal{A}_{\bold{k}_2'}$

4. $g\subseteq f$.

\textbf{Proof: }a. We shall first prove that $\mathcal{A}_{\bold{k}_1} \cap \mathcal{A}_{\bold{k}_2}=\mathcal{A}_{\bold{k}_0}$.
Note that $\mathcal{A}_{\bold{k}_0} \subseteq \mathcal{A}_{\bold{k}_1} \cap \mathcal{A}_{\bold{k}_2}$
is true by the definition of $\leq_2$, so suppose that $\underset{\sim}{a} \in \mathcal{A}_{\bold{k}_1} \setminus \mathcal{A}_{\bold{k}_0}$,
we need to show that $\underset{\sim}{a} \notin \mathcal{A}_{\bold{k}_2}$.
As $\bold{k}_0 \leq_2 \bold{k}_1$, $\underset{\sim}{a}$ is not a
$\mathbb{P}_{\bold{k}_0}$-name. Therefore, it's not a $\mathbb{P}_{\bold{k}_2}$-name,
hence $\underset{\sim}{a} \notin \mathcal{A}_{\bold{k}_2}$.

Now construct $\mathbb P$ as in {[}HwSh:1090{]}, i.e. we take the
amalgamation $\mathbb P'=\mathbb{P}_{\bold{k}_1} \underset{\mathbb{P}_{\bold{k}_0}}{\times} \mathbb{P}_{\bold{k}_2}$
and then we take $\mathbb{P} \in K$ such that $\mathbb{P}' \lessdot \mathbb P$
and $|\mathbb P| \leq (2+|\mathbb{P}'|)^{\aleph_1}$. Now let $\mathcal A:=\mathcal{A}_{\bold{k}_1} \cup \mathcal{A}_{\bold{k}_2}$.
We need to show that $\mathcal A$ is as required, i.e. we need to
prove that $(\mathbb P,\mathcal A)$ satisfy requirements (A)(d) and
(D)(c) in Definition 2 (in the end, we will use claim 5 for the requirement
in Definition (2)(D)(c)). By symmetry, it's enough to show that if
$\Vdash_{\mathbb{P}_{\bold{k}_1}} "\underset{\sim}{b} \in pos(\mathcal{A}_{\bold{k}_1})"$
and $\underset{\sim}{a_0},...,\underset{\sim}{a_{n-1}} \in \mathcal A$,
then $\Vdash_{\mathbb P} "\underset{\sim}{b} \cap (\underset{l<n}{\cap}\underset{\sim}{a_l^{\text{if (l is even)}}}) \in [\omega]^{\omega}"$.
Let $n_2:=n$, wlog there are $n_0<n_1<n_2$ such that $\underset{\sim}{a}_l \in \mathcal{A}_{\bold{k}_0} \iff l<n_0$,
$\underset{\sim}{a}_l \in \mathcal{A}_{\bold{k}_1} \iff l\in [0,n_1)$
and $\underset{\sim}{a}_l \in \mathcal{A}_{\bold{k}_2} \iff l\in [0,n_0) \cup [n_1,n_2)$.
It's enough to show that the last statement is forced by $\mathbb{P}'$,
so let $k<\omega$ and $p=(p_1,p_2) \in \mathbb{P}'$, we shall find
$q\in \mathbb{P}'$ and $m<\omega$ such that $p\leq q$, $k\leq m$
and $q\Vdash_{\mathbb{P}'} " m\in \underset{\sim}{b} \cap (\underset{l<n}{\cap}\underset{\sim}{a_l^{\text{if (l is even)}}})"$.
Let $p_0 \in \mathbb{P}_{\bold{k}_0}$ witness $"(p_1,p_2) \in \mathbb{P}'"$,
i.e. $p_0 \Vdash_{\mathbb{P}_{\bold{k}_0}} "\underset{l=1,2}{\wedge}p_l \in \mathbb{P}_{\bold{k}_l}/\mathbb{P}_{\bold{k}_0}"$.
Let $\underset{\sim}{b^*}=\{m : p_1 \nVdash_{\mathbb{P}_{\bold{k}_1} / \mathbb{P}_{\bold{k}_0}} "m\notin \underset{\sim}{b} \cap (\cap \{\underset{\sim}{a}_l^{\text{if l is even}} : l\in [n_0,n_1)\})"\}$
(so $\underset{\sim}{b^*}$ is a $\mathbb{P}_{\bold{k}_0}$-name)
and let $p_0 \in G_0 \subseteq \mathbb{P}_{\bold{k}_0}$ be generic
over $V$,\textbf{ }then $\underset{\sim}{b^*}=\underset{\sim}{b^*}[G_0] \in V[G_0]$
and as $p_1 \Vdash_{\mathbb{P}_{\bold{k}_1}} "\underset{\sim}{b} \in pos(\mathcal{A}_{\bold{k}_1})"$,
it follows that $p_0 \Vdash_{\mathbb{P}_{\bold{k}_0}} "\underset{\sim}{b^*} \in [\omega]^{\omega}"$,
moreover, $p_0 \Vdash_{\mathbb{P}_{\bold{k}_0}} "\underset{\sim}{b}^* \in pos(\mathcal{A}_{\bold{k}_0})"$.
Let $b^{**}$ be the $\mathbb{P}_{\bold{k}_0}$-name defined as $b^*$
if $p_0$ is in the generic set, and as $\omega$ otherwise. As $\bold{k}_0 \leq_2 \bold{k}_2$,
it follows that $p_2 \Vdash_{\mathbb{P}_{\bold{k}_2}/G_0} "\underset{\sim}{b^{**}} \cap (\underset{l\in [0,n_0) \cup [n_1,n_2)}{\cap}\underset{\sim}{a_l^{\text{if (l is even)}}}) \in [\omega]^{\omega}"$.

Therefore, in $V[G_0]$ there are $(p_2',m)$ such that:

a. $p_2 \leq p_2' \in \mathbb{P}_{\bold{k}_2}/ G_0$.

b. $m>k$.

c. $p_2' \Vdash_{\mathbb{P}_{\bold{k}_2}/ G_0} "m\in \underset{\sim}{b^{**}} \cap (\underset{l \in [0,n_0) \cup [n_1.n_2)}{\cap}\underset{\sim}{a_l^{\text{if (l is even)}}})"$.

Note that as $\underset{\sim}{b^{**}} \in V[G_0]$, $V[G_0] \models "m\in \underset{\sim}{b^{**}}[G_0]=\underset{\sim}{b^*}[G_0]"$.
Therefore, by the definitions of $\underset{\sim}{b}$ and $\underset{\sim}{b^*}$,
there is $p_1' \in \mathbb{P}_{\bold{k}_1}/G_0$ above $p_1$ such
that $p_1' \Vdash_{\mathbb{P}_{\bold{k}_1}/G_0} "m\in \underset{\sim}{b} \cap (\cap \{\underset{\sim}{a}_l^{\text{if l is even}} : l\in [n_0,n_1)\})"$.
Therefore, there is $p_0 \leq p_0' \in G_0$ forcing (in $\mathbb{P}_{\bold{k}_0}$)
all of the aforementioned statements about $(p_1',p_2')$ in $V[G_0]$.
Now it's easy to check that $q=(p_1',p_2')$ is as required. Finally,
extend $(\mathbb P,\mathcal A)$ (with respect to $\leq_1$) to a
member of $K_2^+$. By observation 3, we're done.

b. Follows from (a) by changing names. $\square$

\textbf{Claim 11: }There exists $\bold{k}=(\mathbb{P}_{\bold k},\mathcal{A}_{\bold{k}})=(\mathbb P,\mathcal A) \in K_2^+$
such that $|\mathbb{P}_{\bold k}|=\lambda$ and:

1. For every $X\subseteq \mathbb{P}$ of cardinality $<\mu$, there
exists $\bold{k}'=(\mathbb{Q},\mathcal{A}') \in K_2^+$ such that
$X\subseteq \mathbb Q$, $\bold{k}' \leq_2 \bold k$ and $|\mathbb Q|<\mu$.

2. If $\bold{k}_1,\bold{k}_2 \in K_2^+$, $|\mathbb{P}_{\bold{k}_1}|,|\mathbb{P}_{\bold{k}_2}|<\mu$,
$\bold{k}_1 \leq_2 \bold{k}_2$ and $f_1: \mathbb{P}_{\bold{k}_1} \rightarrow \mathbb P$
is a complete embedding such that $(f_1(\mathbb{P}_{\bold{k}_1}), f_1(\mathcal{A}_{\bold{k}_1}))\leq_2(\mathbb P,\mathcal A)$,
then there is a complete embedding $f_2$ such that $f_1 \subseteq f_2$
and $(f_2(\mathbb{P}_{\bold{k}_2}), f_2(\mathcal{A}_{\bold{k}_2}))\leq_2(\mathbb P,\mathcal A)$. 

\textbf{Proof: }The first property is satisfied by every $\bold k \in K_2^+$
by observation 4. The proof of (2) is as in {[}HwSh:1090{]}. $\square$

\textbf{Claim 12: }A implies B where:

A. a. $\bold{k}_0,\bold{k}_1,\bold{k}_2 \in K_2^+$, $\bold{k}_0 \leq_2 \bold{k}_l$
$(l=1,2)$ and $\mathbb{P}_{\bold{k}_0}=\mathbb{P}_{\bold{k}_1} \cap \mathbb{P}_{\bold{k}_2}$.

b. $\underset{\sim}{D}$ is a $\mathbb{P}_{\bold{k}_0}$-name of a
nonprincipal ultrafilter on $\omega$.

c. For $l=1,2$, $\underset{\sim}{a}_l$ and $\underset{\sim}{b_l}$
are canonical $\mathbb{P}_{\bold{k}_l}$-names of a member of $[\omega]^{\omega}$.

d. For $l=1,2$, $\Vdash_{\mathbb{P}_{\bold{k}_l}} "\underset{\sim}{a}_l \cap \underset{\sim}{b}_l$
is infinite and $\underset{\sim}{a}_l$ contains no members of $\underset{\sim}{D}$
from $V^{\mathbb{P}_{\bold{k}_0}}$.

e. $\mathbb{P}_{\bold{k}_1} \cap \mathbb{P}_{\bold{k}_2}=\mathbb{P}_{\bold{k}_0}$.

f. For $l=1,2$, $\Vdash_{\mathbb{P}_{\bold{k}_l}} "\underset{\sim}{b_l} \text{ is a pseudo intersection of } \underset{\sim}{D}"$.

B. There is $\bold{k} \in K_2^+$ such that:

a. $|\mathbb{P}_{\bold k}| \leq (2+|\mathbb{P}_{\bold{k}_1}|+|\mathbb{P}_{\bold{k}_2}|)^{\aleph_1}$.

b. $\bold{k}_l \leq_2 \bold{k}$ $(l=1,2)$.

c. $\Vdash_{\mathbb{P}_{\bold k}} "\underset{\sim}{a_2} \setminus \underset{\sim}{a_1}$
and $\underset{\sim}{a_1} \setminus \underset{\sim}{a_2}$ are infinite$"$.

\textbf{Proof: }Let $\bold{k} \in K_2^+$ be the object constructed
by the proof of claim 10. We need to prove that $\bold k$ satisfies
clause (B)(c). For $l=0,1,2$, let $\mathbb{P}_l=\mathbb{P}_{\bold{k}_l}$
and let $\mathbb P'$ be as in the proof of claim 10, so it suffices
to prove that $\Vdash_{\mathbb{P}'} "\underset{\sim}{a_2} \setminus \underset{\sim}{a_1}$
and $\underset{\sim}{a_1} \setminus \underset{\sim}{a_2}$ are infinite$"$.
Let $p=(p_1,p_2) \in \mathbb{P}'$, $k<\omega$ and let $p_0 \in \mathbb{P}_0$
be a witness of $"p=(p_1,p_2) \in \mathbb{P}'"$. Now let $G_0 \subseteq \mathbb{P}'$
be generic over $V$ such that $p_0 \in G_0$ and let $D=\underset{\sim}{D}[G_0]$.
By the assumptions, for $l=1,2$, $p_l \Vdash_{\mathbb{P}_1 /G_0} "\underset{\sim}{a_l} \cap \underset{\sim}{b_l}$
is an infinite pseudo intersection of $D"$. In $V[G_0]$, let $b_l^*=\{m: p_l \nVdash_{\mathbb{P}_l/G_0} "m\notin \underset{\sim}{a_l} \cap \underset{\sim}{b_l}"\}$,
then $p_l \Vdash_{\mathbb{P}_l/G_0} "\underset{\sim}{a_l} \cap \underset{\sim}{b_l} \subseteq \underset{\sim}{b_l^*}$,
hence $b_l^*$ infinite$"$. As $p_l \Vdash_{\mathbb{P}_l /G_0} "\underset{\sim}{a_l} \cap \underset{\sim}{b_l}$
is a pseudo intersection of $D"$, necessarily $V[G_0] \models b_l^* \in D$.

Let $a_l^*=\{m : p_l \Vdash_{\mathbb{P}_l/ G_0} m\in \underset{\sim}{a_l}\}$,
so $a_l^* \in V[G_0]$ and $p_l \Vdash_{\mathbb{P}_l/G_0} "a_l^* \subseteq \underset{\sim}{a_l}"$.
Now recall that $\Vdash_{\mathbb{P}_{\bold{k}_l}} "\underset{\sim}{a}_l$
contains no member of $\underset{\sim}{D}$ from $V^{\mathbb{P}_{\bold{k}_0}}"$,
therefore $p_l \Vdash_{\mathbb{P}_l /G_0} "a_l^* \notin \underset{\sim}{D}"$. 

Hence in $V[G_0]$ (recalling $D$ is an ultrafilter), $b:=(b_1^* \cap b_2^*) \setminus (a_1^* \cup a_2^*) \in D$.
Let $m\in b$ be such that $k<m$. By the definition of $b_l^*$,
there is $p_l' \in \mathbb{P}_l /G_0$ above $p_l$ such that $p_l' \Vdash_{\mathbb{P}_l /G_0} "m\in \underset{\sim}{a_l} \cap \underset{\sim}{b_l}"$.
By the definition of $a_l^*$, there is $p_l'' \in \mathbb{P}_l /G_0$
above $p_l$ such that $p_l'' \Vdash_{\mathbb{P}_l /G_0} "m\notin \underset{\sim}{a_l}"$.
Let $p_0' \in G_0$ be a condition above $p_0$ forcing the above
statements, so $p_0'$ is witnessing the fact that $(p_1',p_2''),(p_1'',p_2') \in \mathbb{P}'$
are above $p=(p_1,p_2)$. Now $m>k$, $(p_1',p_2'') \Vdash " m\in \underset{\sim}{a_1} \setminus \underset{\sim}{a_2}"$
and $(p_1'',p_2') \Vdash "m\in \underset{\sim}{a_2} \setminus \underset{\sim}{a_1}"$,
which completes the proof. $\square$

\textbf{Definition 13: }Let $\mathbb P=\mathbb{P}_{\bold k}$ be the
forcing from claim 11, let $G\subseteq \mathbb P$ be generic over
$V$ and in $V[G]$, let $V_1=HOD(\mathbb{R}^{<\mu} \cup \{\mathcal{A}_{\bold k}\})$.

\textbf{Definition 14 ({[}NN{]}): }A family $\mathcal F \subseteq [\omega]^{\omega}$
is nearly mad if $|A\cap B|<\aleph_0$ or $|A\Delta B|<\aleph_0$
for every $A\neq B \in \mathcal F$, and $\mathcal F$ is maximal
with respect to this property.

\textbf{Theorem 15: }$V_1 \models ZF+DC_{<\mu}+"\text{there are no mad families}"+"\text{there are no nearly mad families}"+"\text{there exists a non-meager filter on } \omega"$.

\textbf{Proof: }1. In order to see that there exists a non-meagre
filter in $V_1$, let $\underset{\sim}{D}$ be the filter generated
by $\mathcal{A}_{\bold k}$ and the cofinite sets. By claim 8 and
the choice of $\bold k$, $\underset{\sim}{D}$ contains many Cohen
reals and therefore is non-meager.

2. The proof of the non-existence of mad families is exactly as in
{[}HwSh:1090{]}, where $(K_2^+,\leq_2)$ here replaces $(K,\lessdot)$
there, and claim 10 is used for the amalgamation arguments. Alternatively,
see the proof of (3) below.

3. The non-existence of a nearly mad family in $V_1$ will follow
from the proofs below. $\square$

\noindent \begin{center}
\textbf{Somewhat mad families}
\par\end{center}

\textbf{Definition 16: }A family $\mathcal F \subseteq [\omega]^{\omega}$
is somewhat mad if:

a. For every $a_1,a_2 \in \mathcal F$, $|a_1 \cap a_2|<\aleph_0$
or $a_1 \subseteq^* a_2$ or $a_2 \subseteq^* a_1$.

b. If $b\in [\omega]^{\omega}$ then for some $a\in \mathcal F$,
$|a\cap b|=\aleph_0$.

\textbf{Observation 17: }Nearly mad families are somewhat mad. $\square$

\textbf{Definition 18: }Let $Pr(\bold{k}_1,\bold{k}_2,\underset{\sim}{D},\underset{\sim}{b_2})$
mean:

a. $\bold{k}_1 \leq_2 \bold{k}_2$, $\underset{\sim}{b_2}$ is a $\mathbb{P}_{\bold{k}_2}$-name
and $\underset{\sim}{D}$ is a $\mathbb{P}_{\bold{k}_1}$-name such
that $\Vdash_{\mathbb{P}_{\bold{k}_2}} "\underset{\sim}{b_2} \in [\omega]^{\omega}"$
and $\Vdash_{\mathbb{P}_{\bold{k}_1}} "\underset{\sim}{D}$ is a nonprincipal
ultrafilter on $\omega"$.

b. If $G_1 \subseteq \mathbb{P}_{\bold{k}_1}$ is generic over $V$,
$p_1 \in \mathbb{P}_{\bold{k}_2}/G_1$, $b_0^*=\{n : p_1 \Vdash_{\mathbb{P}_{\bold{k}_1}G_1} "n\in \underset{\sim}{b_2}" \}$
and $b_1^*=\{n : p_1 \nVdash_{\mathbb{P}_{\bold{k}_1}/G_1} "n\notin \underset{\sim}{b_2}" \}$,
then $V[G_1] \models b_1^* \setminus b_0^* \in D$.

\textbf{Claim 19: }(A) implies (B) where:

A. a. $\bold{k}_1 \in K_2^+$.

b. $\Vdash_{\mathbb{P}_{\bold{k}_1}} "\underset{\sim}{D}$ is a nonprincipal
ultrafilter on $\omega"$.

c. $\Vdash_{\mathbb{P}_{\bold{k}_1}} "\underset{\sim}{S_1}$ is somewhat
mad$"$.

d. $\Vdash_{\mathbb{P}_{\bold{k}_1}} " \underset{\sim}{S_1} \cap \underset{\sim}{D}=\emptyset"$.

B. There is $\bold{k}_2$ such that:

a. $\bold{k}_2 \in K_2^+$.

b. $\bold{k}_1 \leq_2 \bold{k}_2$.

c. $|\mathbb{P}_{\bold{k}_2}| \leq (2+|\mathbb{P}_{\bold{k}_1}|)^{\aleph_1}$.

d. ($\alpha$) implies ($\beta$) where:

$\alpha$. $(\bold{k}_3,\underset{\sim}{S_2})$ satisfy the following
properties:

1. $\bold{k}_2 \leq_2 \bold{k}_3 \in K_2^+$.

2. $\Vdash_{\mathbb{P}_{\bold{k}_3}} "\underset{\sim}{S_2}$ is somewhat
mad and $\underset{\sim}{S_1} \subseteq \underset{\sim}{S_2}"$.

3. $\Vdash_{\mathbb{P}_{\bold{k}_3}} "$no member of $\underset{\sim}{S_2} \setminus \underset{\sim}{S_1}$
contains a member of $\underset{\sim}{D}"$.

$\beta$. For some $\mathbb{P}_{\bold{k}_3}$-name $\underset{\sim}{a}$,
$\Vdash_{\mathbb{P}_{\bold{k}_3}} "\underset{\sim}{a} \in \underset{\sim}{S_2}"$
and $Pr(\bold{k}_1,\bold{k}_3,\underset{\sim}{D},\underset{\sim}{a})$.

\textbf{Proof: }Using Mathias forcing restricted to $\underset{\sim}{D}$,
it's easy to see that there is $\bold{k}_2$ and a $\mathbb{P}_{\bold{k}_2}$
name $\underset{\sim}{b}$ such that $\bold{k}_1 \leq_2 \bold{k}_2 \in K_2^+$,
$|\mathbb{P}_{\bold{k}_2}| \leq (2+|\mathbb{P}_{\bold{k}_1}|)^{\aleph_1}$
and $\Vdash_{\mathbb{P}_{\bold{k}_2}} "\underset{\sim}{b}$ is a pseudo
intersection of $\underset{\sim}{D}"$. Therefore, $\Vdash_{\mathbb{P}_{\bold{k}_2}} "\underset{\sim}{b} \in [\omega]^{\omega}$
is almost disjoint to every $\underset{\sim}{a} \in \underset{\sim}{S_1}"$
(by the fact that $\Vdash_{\mathbb{P}_{\bold{k}_1}} "\underset{\sim}{S_1} \cap \underset{\sim}{D}=\emptyset"$).

We shall now prove that $\bold{k}_2$ satisfies (B)(d). Suppose that
$(\bold{k}_3,\underset{\sim}{S_2})$ are as there. By the somewhat
madness of $\underset{\sim}{S_2}$, there is $\underset{\sim}{a}$
such that $\Vdash_{\mathbb{P}_{\bold{k}_3}} "\underset{\sim}{a} \in \underset{\sim}{S_2}$
and $|\underset{\sim}{a} \cap \underset{\sim}{b}|=\aleph_0"$. Therefore,
$\Vdash_{\mathbb{P}_{\bold{k}_3}} "\underset{\sim}{a} \cap \underset{\sim}{b} \in [\omega]^{\omega}$
is a pseudo intersection of $\underset{\sim}{D}"$. Now let $G_1 \subseteq \mathbb{P}_{\bold{k}_1}$
be generic over $V$. If $p_1 \in \mathbb{P}_{\bold{k}_3} /G_1$ then
$b^*=\{n: p_1 \nVdash_{\mathbb{P}_{\bold{k}_3}/G_1} "n\notin \underset{\sim}{a} \cap \underset{\sim}{b}" \} \in \underset{\sim}{D}[G_1]$
by the fact that $\underset{\sim}{a} \cap \underset{\sim}{b}$ is
a pseudo intersection of $\underset{\sim}{D}$. Now let $a^*=\{n: p_1 \Vdash_{\mathbb{P}_{\bold{k}_3 / G_1}} "n\in \underset{\sim}{a}"\}$,
then $p_1 \Vdash_{\mathbb{P}_{\bold{K}_3} /G_1} "a^* \subseteq \underset{\sim}{a} \text{ is infinite}"$.
If $a^* \in \underset{\sim}{D}[G_1]$, then $p_1$ forces that $\underset{\sim}{a}$
(which belongs to $\underset{\sim}{S_2} \setminus \underset{\sim}{S_1}$)
contains a member of $\underset{\sim}{D}[G_1]$, contradicting $(\alpha)(3)$.
Therefore, $a^* \notin \underset{\sim}{D}[G_1]$, and $\underset{\sim}{a}$
is as required in the definition of $Pr(\bold{k}_1,\bold{k}_3,\underset{\sim}{D},\underset{\sim}{a})$.
$\square$

\textbf{Claim 20: }There is no somewhat mad family in $V_1$.

\textbf{Proof: }Suppose towards contradiction that $\underset{\sim}{S}$
is a $\mathbb{P}$-name of a somewhat mad family. As in {[}HwSh:1090{]},
let $\underset{\sim}{D}$ be a $\mathbb P$-name of a Ramsey ultrafilter
on $\omega$ such that $\Vdash_{\mathbb P} "\underset{\sim}{S} \cap \underset{\sim}{D}=\emptyset"$.
By claim 11(a), there is $\bold{k}_1 \leq_2 \bold k$ such that $\bold{k}_1 \in K_2^+$,
$|\mathbb{P}_{\bold{k}_1}|<\mu$ and $\underset{\sim}{S}$ is definable
using a $\mathbb{P}_{\bold{k}_1}-$name. Let $K_{\mathbb P}^+$ be
the set of $\bold{k}' \in K_2^+$ such that $\bold{k}' \leq_2 \bold k$,
$|\mathbb{P}_{\bold{k}'}|<\mu$, $\underset{\sim}{S} \restriction \mathbb{P}_{\bold{k}'}$
is a canonical $\mathbb{P}_{\bold{k}'}$-name of a somewhat mad family
in $V^{\mathbb{P}_{\bold{k}'}}$ and $\underset{\sim}{D} \restriction \mathbb{P}_{\bold{k}'}$
is a $\mathbb{P}_{\bold{k}'}$-name of a Ramsey ultrafilter on $\omega$.
As in {[}HwSh:1090{]}, $K_{\mathbb P}^+$ is $\leq_2$-dense in $K_2^+$,
so there exists $\bold{k}_2 \in K_{\mathbb P}^+$ such that $\bold{k}_1 \leq \bold{k}_2$.
Let $\bold{k}_3 \in K_2^+$ be as in claim 19 for $(\bold{k}_2,\underset{\sim}{S} \restriction \mathbb{P}_{\bold{k}_2})$,
wlog $\bold{k}_3 \leq_2 \bold k$ (see claim 11). Choose $\bold{k}_4 \in K_{\mathbb P}^+$
such that $\bold{k}_3 \leq_2 \bold{k}_4$ and let $\underset{\sim}{S_4}:=\underset{\sim}{S} \restriction \mathbb{P}_{\bold{k}_4}$.
Let $\underset{\sim}{a}$ be a $\mathbb{P}_{\bold{k}_4}$-name such
that $Pr(\bold{k}_3,\bold{k}_4,\underset{\sim}{D},\underset{\sim}{a})$
holds, as guaranteed by claim 19.

As in {[}HwSh:1090{]}, there are $\bold{k}_5,\bold{k}_6, \in K_{\mathbb P}^+$
and an isomorphism $f$ from $\bold{k}_4$ to $\bold{k}_5$ over $\bold{k}_2$
such that $\mathbb{P}_{\bold{k}_5}$ adds a generic for $\mathbb{M}_{\underset{\sim}{D} \restriction \mathbb{P}_{\bold{k}_4}}$
(Mathias forcing restricted to the ultrafilter $\underset{\sim}{D} \restriction \mathbb{P}_{\bold{k}_4}$)
and $(\bold{k}_2,\bold{k}_4,\bold{k}_5,\bold{k}_6)$ here are as $(\bold{k}_0,\bold{k}_1,\bold{k}_2,\bold{k}_3)$
in claim 10, and wlog $\bold{k}_6 \leq_2 \bold k$. By the choice
of $f$, $\Vdash_{\mathbb P} "\underset{\sim}{a},f(\underset{\sim}{a}) \in \underset{\sim}{S}"$.

By claim 12, with $(\bold{k}_2,\bold{k}_4,\bold{k}_5,\bold{k}_6,\underset{\sim}{a},f(\underset{\sim}{a}))$
standing for $(\bold{k}_0,\bold{k}_1,\bold{k}_2,\bold{k}_,\underset{\sim}{a}_1,\underset{\sim}{a}_2)$
there, it's forced by $\mathbb{P}_{\bold{k}_5}$, and hence by $\mathbb P$,
that $\underset{\sim}{a} \setminus f(\underset{\sim}{a})$ and $f(\underset{\sim}{a}) \setminus \underset{\sim}{a}$
are infinite. As in {[}HwSh:1090{]}, $\Vdash_{\mathbb P} "|\underset{\sim}{a} \cap f(\underset{\sim}{a})|=\aleph_0"$.
As $\Vdash_{\mathbb P} "\underset{\sim}{a},f(\underset{\sim}{a}) \in \underset{\sim}{S}"$,
we get a contradiction. $\square$

\textbf{\large References}{\large \par}

{[}HwSh:1090{]} Haim Horowitz and Saharon Shelah, Can you take Toernquist's
inaccessible away? arXiv:1605.02419

{[}Ma1{]} A. R. D. Mathias, Happy families, Ann. Math. Logic \textbf{12
}(1977), no. 1, 59-111.

{[}Ma2{]} A. R. D.. Mathias, A remark on rare filters.

{[}NN{]} Itay Neeman and Zach Norwood, Happy and mad families in $L(\mathbb R)$,
preprint.

{[}To{]} Asger Toernquist, Definability and almost disjoint families,
arXiv:1503.07577.

$\\$

(Haim Horowitz) Einstein Institute of Mathematics

Edmond J. Safra campus,

The Hebrew University of Jerusalem.

Givat Ram, Jerusalem, 91904, Israel.

E-mail address: haim.horowitz@mail.huji.ac.il

$\\$

(Saharon Shelah) Einstein Institute of Mathematics

Edmond J. Safra campus,

The Hebrew University of Jerusalem.

Givat Ram, Jerusalem, 91904, Israel.

Department of Mathematics

Hill Center - Busch Campus,

Rutgers, The State University of New Jersey.

110 Frelinghuysen road, Piscataway, NJ 08854-8019 USA

E-mail address: shelah@math.huji.ac.il
\end{document}